\date{}
\documentclass [11pt]{article}
\usepackage{amsfonts,amsmath,amssymb,amsthm,graphicx}

\title{The Tur\'an number of $F_{3,3}$}
\author{
Peter Keevash \thanks{School of Mathematical Sciences,
Queen Mary, University of London, Mile End Road, London E1 4NS, UK.
Email: p.keevash@qmul.ac.uk.
Research supported in part by
ERC grant 239696 and
EPSRC grant EP/G056730/1.}
\and Dhruv Mubayi \thanks{Department of Mathematics, Statistics, and Computer
Science, University of Illinois,  Chicago, IL 60607. Email: mubayi@math.uic.edu. Research
supported in part by NSF grant DMS-0969092. }
}

\theoremstyle{plain}
\newtheorem{theo}{Theorem}[section] 

\newtheorem{lemma}[theo]{Lemma}

\theoremstyle{definition}

\newcommand{\nib}[1]{\noindent {\bf #1}}

\newcommand{\sub}{\subseteq}

\newcommand{\ex}{\mbox{ex}}

\newcommand{\sm}{\setminus}

\def\COMMENT#1{}

\def\qed{\hfill $\Box$}

\begin{document}
\maketitle

\begin{abstract}
Let $F_{3,3}$ be the $3$-graph on $6$ vertices,
labelled $abcxyz$, and $10$ edges, one of which is $abc$,
and the other $9$ of which are all triples that
contain $1$ vertex from $abc$ and $2$ vertices from $xyz$.
We show that for all $n \ge 6$, the maximum number of edges
in an $F_{3,3}$-free $3$-graph on $n$ vertices is
$\binom{n}{3} -  \binom{\lfloor n/2 \rfloor}{3}
- \binom{\lceil n/2 \rceil}{3}$.
This sharpens results of Zhou \cite{Z91} and of
the second author and R\"odl \cite{MR02}.
\end{abstract}

\section{Introduction}

The {\em Tur\'an number} $\ex(n,F)$ is the maximum number of edges
in an $F$-free $r$-graph on $n$ vertices.%
\footnote{
An {\em $r$-graph} (or {\em $r$-uniform hypergraph}) $G$ consists of a vertex set
and an edge set, each edge being some $r$-set of vertices.
We say $G$ is {\em $F$-free} if it does not have a (not necessarily induced)
subgraph isomorphic to $F$.
}
It is a long-standing open problem in Extremal Combinatorics
to understand these numbers for general $r$-graphs $F$.
For ordinary graphs ($r=2$) the picture is fairly complete,
although there are still many open problems, such as determining
the order of magnitude for Tur\'an numbers of bipartite graphs.
However, for $r \ge 3$ there are very few known results.
Having solved the problem for the complete graph $F=K_t$,
Tur\'an \cite{T61} posed the natural question of determining
$\ex(n,F)$ when $F=K^r_t$ is a complete $r$-graph on $t$ vertices.
To date, no case with $t>r>2$ of this question has been solved, even asymptotically.
Despite the lack of progress on the Tur\'an problem for complete hypergraphs,
there are certain hypergraphs for which the problem has been solved asymptotically,
or even exactly; we refer the reader to the survey \cite{K11}. While it would be more
satisfactory to have a general theory, we may hope that this will develop out
of the methods discovered in solving isolated examples.

The contribution of this paper is a surprisingly short complete solution
to the Tur\'an problem for the following $3$-graph.
Let $F_{3,3}$ be the $3$-graph on $6$ vertices,
labelled $abcxyz$, and $10$ edges, one of which is $abc$,
and the other $9$ of which are all triples that
contain $1$ vertex from $abc$ and $2$ vertices from $xyz$.
A lower bound for $\ex(n,F_{3,3})$ is given by the following natural construction.
Let $B(n)$ denote the balanced complete bipartite $3$-graph,
which is obtained by partitioning a set of $n$ vertices into
parts of size $\lfloor n/2 \rfloor$ and $\lceil n/2 \rceil$,
and taking as edges all triples that are not contained within
either part. Let
\[ b(n) = \binom{n}{3} -  \binom{\lfloor n/2 \rfloor}{3}
- \binom{\lceil n/2 \rceil}{3}\]
denote the number of edges in $B(n)$.
We prove the following result.

\begin{theo} \label{main}
For any $n \ge 1$, $\ex(n,F_{3,3}) = b(n)$, unless $n=5$, when $\ex(5,F_{3,3}) = 10$.
\end{theo}

The Tur\'an problem for $F_{3,3}$ was previously studied by Mubayi and R\"odl \cite{MR02},
who obtained the asymptotic result $\ex(n,F_{3,3}) = (1+o(1))b(n)$.
Another related problem is the following result of Zhou \cite{Z91}.
Say that two vertices $x,y$ in a $3$-graph $G$ are
{\em $t$-connected} if there are vertices $a,b,c$ such that every triple
with $2$ vertices from $abc$ and $1$ from $xy$ is an edge. Say that $xyz$
is a {\em $t$-triple} if $xyz$ is an edge and each pair in $xyz$ is $t$-connected.
For example, $K^3_5$ is a $t$-triple (this is the motivation for the definition).
The result of  \cite{Z91} is that the unique largest
$3$-graph on $n$ vertices with no $t$-triple is complete bipartite.
Note that $F_{3,3}$ is a $t$-triple, so Theorem \ref{main} strengthens
Zhou's extremal result (but not the classification of the extremal example).

Our proof uses the link multigraph method introduced by de Caen and F\"uredi \cite{dCF00}.
There are now a few examples where this method has been used to obtain asymptotic
results, or exact results for $n$ sufficiently large. We used it in \cite{KM04} to obtain
an exact result for cancellative $3$-graphs for all $n$, and an exact result
for the configuration $F_5 = \{123,124,345\}$ for $n \ge 33$.
However, until recently there were no known applications to an
exact result for all $n$ with a single forbidden hypergraph.
The result in this paper gives such an application;
another was given very recently by Goldwasser \cite{G+},
who obtained an exact result for $F_5$ for all $n$.

In the next section we describe the link multigraph construction
and prove a lemma that applies to such multigraphs. We use this
to prove Theorem \ref{main} in Section 3. The final section
contains some concluding remarks and open problems.

\section{A multigraph lemma}

The proof of Theorem \ref{main} will use the following construction of a multigraph
from a $3$-graph $G$. Suppose $S$ is a set of vertices in $G$.
The `link multigraph' of $S$ has vertex set $X = V \sm S$ and edge set
$M = \sum_{a \in S} G(a)[X]$. Here we write $G(a) = \{xy: axy \in G\}$,
denote the restriction to $X$ by $[X]$, and use summation to denote multiset union.
Thus we obtain a multigraph $M$ in which each pair of vertices has multiplicity between $0$ and $|S|$.
Furthermore, we may regard each edge of $M$ as being `coloured' by a vertex in $S$
(an edge may have several colours). We write $w(xy)$ for the multiplicity of the pair $xy$ in $M$.

Now suppose $M$ is any multigraph on $n$ vertices (not necessarily as above).
Write $w(xy)$ for the multiplicity of a pair $xy$ in $M$,
and write $e(M)$ for the sum of $w(xy)$ over all
(unordered) pairs of vertices in $M$. For any $S \sub V(M)$
let $i(S)$ denote the sum of $w(xy)$ over all
pairs of vertices that contain at least one vertex of $S$.
If $S=\{x\}$ consists of a single vertex then $i(S)=d(x)$
is the weighted degree of $x$.
Define
\begin{equation*}
m(n) = \begin{cases} \frac{3n^2}{2}-n & \text{if $n$ is even,} \\
\frac{3n^2-1}{2}-n & \text{if $n$ is odd.} \end{cases}
\end{equation*}

\begin{lemma}\label{multi}
Suppose $M$ is a multigraph on $n$ vertices
with $0 \le w(xy) \le 4$ for every pair $xy$
and $w(xy)+w(xz)+w(yz) \le 10$ for every triple $xyz$.
Then $e(M) \le m(n)$.
\end{lemma}

\nib{Proof.}
We argue by induction on $n$. The statement is trivial for $n=1$
and $n=2$, and is immediate from the assumption on triples for $n=3$.
Now suppose that $n \ge 4$. We consider separate cases according to the
parity of $n$.

Suppose first that $n$ is even, $M$ is a multigraph satisfying the
hypotheses of the lemma, and suppose for a contradiction that
$e(M) = \frac{3n^2}{2}-n+1$. Since $e(M) > 3\binom{n}{2}$ we
can choose a pair $xy$ with $w(xy)=4$. If we delete $xy$ then
we obtain a multigraph $M'$ on $n-2$ satisfying the hypotheses of the lemma.
By induction hypothesis we have $e(M') \le \frac{3(n-2)^2}{2}-(n-2)$.
Then $i(xy) = e(M)-e(M') \ge \frac{3n^2}{2}-n+1 - \frac{3(n-2)^2}{2} + (n-2) = 6n-7$.
Now the sum of $w(xz)+w(yz)$ over $z$ in $V(M) \sm \{x,y\}$ is
$i(xy)-w(xy) \ge 6n-11 > 6(n-2)$, so there must be some $z$
with $w(xz)+w(yz) \ge 7$. But then $w(xy)+w(xz)+w(yz) \ge 11$, contradiction.

The argument for $n$ odd is similar. Suppose for a contradiction that $M$
satisfies the hypotheses of the lemma but $e(M) = \frac{3n^2-1}{2}-n+1$.
Choose $xy$ with $w(xy)=4$. The induction hypothesis gives
$i(xy) \ge \frac{3n^2-1}{2}-n+1 - \frac{3(n-2)^2-1}{2} + (n-2) = 6n-7$.
Then, as in the case of $n$ even, we have $i(xy)-w(xy) > 6(n-2)$,
so there must be some $z$ with $w(xz)+w(yz) \ge 7$, contradiction. \qed

\medskip

Note the following two examples where equality holds in Lemma \ref{multi}.
(We do not claim that these are the only cases of equality.)
\begin{enumerate}
\item
Define a multigraph $M_1(n)$ on $n$ vertices as follows.
Let $A \cup B$ be a balanced partition of the vertex set.
Let crossing pairs have multiplicity $4$
and pairs inside each part have multiplicity $2$.
If $n$ is even then $e(M_1(n)) = 2\binom{n}{2} + 2(n/2)^2 = 3n^2/2-n$.
If $n$ is odd then $e(M_1(n)) = 2\binom{n}{2} + 2 \frac{n^2-1}{4} = \frac{3n^2-1}{2}-n$.
\item
Define a multigraph $M_2(n)$ on $n$ vertices as follows.
Let all pairs have multiplicity $3$ except for a maximum
size matching of multiplicity $4$.
If $n$ is even then $e(M_2(n)) = 3\binom{n}{2} + n/2 = 3n^2/2-n$.
If $n$ is odd then $e(M_2(n)) = 3\binom{n}{2} + \frac{n-1}{2} = \frac{3n^2-1}{2}-n$.
\end{enumerate}

The following two calculations will also be useful.
Note that $M_2(n-1)$ can be obtained from $M_2(n)$ by deleting a vertex,
which can be any vertex if $n$ is even, but must be the vertex not
incident to an edge of multiplicity $4$ when $n$ is odd.
Then $m(n)-m(n-1)$ is equal to the number of edges removed,
which is $3(n-2)+4 = 3(n-1)+1$ when $n$ is even,
or $3(n-1)$ when $n$ is odd.

Next consider a copy of $B(n)$ with parts $A$ and $B$.
Construct a copy of $B(n-4)$ by removing vertices $wx$ from $A$ and $yz$ from $B$.
Then we have $b(n)-b(n-4) = i(wxyz)$, where similarly to our multigraph notation,
we let $i(S)$ denote the number of edges that contain at least one vertex of $S$.
We can count $i(wxyz)$ as follows.  There are $4$ edges of $B(n)$ contained in $wxyz$.
Next consider edges with $2$ vertices in $wxyz$. The $4$ crossing pairs $wy$, $wz$,
$xy$, $xz$ form an edge with each of the $n-4$ vertices of $V \sm \{w,x,y,z\}$,
so contribute $4(n-4)$ edges. The pairs $wx$ and $yz$ form an edge with any of the
vertices in the other part, so contribute $n-4$ edges (this holds whether $n$ is even or odd).
Finally, note that the link multigraph of $wxyz$ in $B(n)$ is precisely $M_1(n-4)$.
Thus the number of edges with $1$ vertex in $wxyz$ is $m(n-4)$.
Then we have
\[b(n)-b(n-4) = m(n-4) + 5(n-4) + 4.\]

\section{Proof of Theorem \ref{main}}

We start with the lower bound. We have already described the construction $B(n)$,
but we also need to check that it is $F_{3,3}$-free. To see this, suppose for a
contradiction that there is a copy of $F_{3,3}$ in $B(n)$, labelled $abcxyz$ as above.
Label the parts of $B(n)$ as $X$ and $Y$, and suppose without loss of generality
that $a \in X$. The edges $axy$, $axz$, $ayz$ can only be simultaneously realised
by putting all of $x,y,z$ in $Y$, or $2$ of $x,y,z$ in $Y$ and $1$ in $X$.
Either way, the edges $bxy$, $bxz$, $byz$ imply that $b$ is in $X$,
and the edges $cxy$, $cxz$, $cyz$ imply that $c$ is in $X$.
But this contradicts the fact that $abc$ is an edge. Thus $B(n)$ is $F_{3,3}$-free.
This shows that $\ex(n,F_{3,3}) \ge b(n)$. This bound can be improved when $n=5$,
as $b(5)=9$, but the complete $3$-graph $K^3_5$ is $F_{3,3}$-free.
It follows that $\ex(5,F_{3,3}) = 10$ (obviously it cannot be larger).

The main task in the proof is to establish the upper bound.
We prove the following statement by induction on $n$:
\begin{itemize}
\item
Suppose $G$ is an $F_{3,3}$-free $3$-graph on $n \ge 1$ vertices.
Then $e(G) \le b(n)$, unless $n=5$, in which case $e(G) \le 10$.
\end{itemize}
Note that this statement is trivial for $n \le 5$,
as $B(n)$ is complete for $n \le 4$, and for $n=5$
the statement allows $e(K^3_5)=10$ edges.
Furthermore, the bound holds for $n=6$, as $B(6)$ is
obtained by deleting $2$ edges from $K^3_6$,
and it is clear that if one only deletes one edge from $K^3_6$
then there is a copy of $F_{3,3}$. Moreover, $B(6)$ is the unique
$F_{3,3}$-free $3$-graph on $6$ vertices with $\binom{6}{3}-2 = 18$ edges.
To see this, we exhibit appropriate edges in the complement of $F_{3,3}$ as defined
above: $xab$ and $yac$ are not in $F_{3,3}$ and intersect in $1$ vertex,
whereas  $xab$ and $xac$ are not in $F_{3,3}$ and intersect in $2$ vertices.

Now suppose for a contradiction that  $n \ge 7$ with $n \ne 9$ and
$G$ is an $F_{3,3}$-free $3$-graph on $n$ vertices with $e(G) = b(n)+1$.
(We will need to modify the argument in the case $n=9$.)
We start by finding a copy of $K^3_4$ in $G$.
For this we use the following averaging argument.
Given a $3$-graph $H$ on $m$ vertices, let $d(H) = e(H)\binom{m}{3}^{-1}$
denote the density of $H$. A simple calculation shows that $d(G)$ is
the average of $d(G \sm v)$ over all vertices $v$ of $G$. Note that
deleting a vertex from $B(n)$ leaves a complete bipartite $3$-graph on $n-1$
vertices; it is not necessarily balanced, but certainly has at most $b(n-1)$ edges.
It follows that $d(B(n)) \le d(B(n-1))$, i.e.\ $d(B(n))$ is non-increasing in $n$.
Since $d(B(n)) \to 3/4$ as $n \to \infty$ we have $d(B(n)) \ge 3/4$ for all $n$.
Since $e(G) > b(n)$ we have $d(G) > 3/4$. Averaging again, we see that there
is a set $abcd$ of $4$ vertices where $d(G[abcd])>3/4$. This implies that all
$4$ triples in $abcd$ are edges of $G$, as desired.

Note that $G \sm \{a,b,c,d\}$ is
an $F_{3,3}$-free $3$-graph on $n-4$ vertices with $e(G) - i(abcd)$ edges.
By induction this is at most $b(n-4)$ (since $n \ne 9$),
so we obtain $i(abcd) \ge b(n)-b(n-4)+1$.
Now we count the edges incident to $abcd$ according to the number of vertices of $abcd$
they contain. There are $4$ such edges contained in $abcd$.
To estimate edges with one vertex in $abcd$ let $M$ be the link multigraph of $abcd$ in $G$.
Note that there is no triangle $xyz$ in $M$ such that each pair $xy$, $xz$, $yz$
is coloured by the same set of $3$ colours from $abcd$: this would give a copy of $F_{3,3}$.
This implies that $w(xy)+w(xz)+w(yz) \le 10$ for every triple $xyz$ in $M$.
Thus we can apply Lemma \ref{multi} to get $e(M) \le m(n-4)$.
We conclude that the number of edges with $2$ vertices in $abcd$
is at least $b(n)-b(n-4)+1 - 4 - m(n-4) = 5(n-4)+1$.
It follows that there is some $e \in V(G) \sm \{a,b,c,d\}$
such that all $6$ pairs from $abcd$ form an edge with $e$.
Thus $abcde$ forms a copy of $K^3_5$ in $G$.

For each $x \in abcde$ we have $i(abcde \sm x) \ge b(n)-b(n-4)+1$, so
\[\Sigma := \sum_{x \in abcde} i(abcde \sm x) \ge 5(b(n)-b(n-4)+1) = 5(m(n-4)+5(n-3)).\]
We can also count $\Sigma$ according to the intersection of edges with $abcde$.
Edges with at least $2$ vertices in $abcde$ are counted $5$ times,
and edges with $1$ vertex in $abcde$ are counted $4$ times.
By Lemma \ref{multi}, for each $x \in abcde$ the link multigraph of $abcde \sm x$
restricted to $V(G) \sm \{a,b,c,d,e\}$ has at most $m(n-5)$ edges. By averaging, there
are at most $\frac{5}{4}m(n-5)$ edges with $1$ vertex in $abcde$.
Since these are counted $4$ times they contribute at most $5m(n-5)$ to $\Sigma$.
Also, $abcde$ is complete, so we have $10$ edges inside $abcde$.

Writing $Z$ for the number of edges with $2$ vertices in $abcde$, we obtain
\[ 5(m(n-4)+5(n-3)) = 5(b(n)-b(n-4)+1) \le \Sigma \le 5(10+Z+m(n-5)),\]
so $Z \ge m(n-4)-m(n-5) + 5(n-5)$.
Recall that $m(n-4)-m(n-5)$ is $3(n-5)+1$ when $n$ is even,
or $3(n-5)$ when $n$ is odd. Thus $Z \ge 8(n-5)$.
It follows that there is some $f \in V(G) \sm \{a,b,c,d,e\}$
such that at least $8$ pairs from $abcde$ form an edge with $f$.
Thus $abcdef$ is obtained from $K^3_6$ by deleting at most $2$ edges,
and if $2$ edges are deleted they cannot be disjoint, as they both contain $f$.
As noted above, this implies that $abcdef$ contains $F_{3,3}$, so we have a contradiction.

It remains to prove the bound for $n=9$. Again we start by choosing
$abcd$ as a copy of $K^3_4$ in $G$. Since $b(5)=\binom{5}{3}-1$,
we obtain $i(abcd) \ge b(n)-b(n-4)$. Then the same calculation as above
shows that there are at least $5(n-4)$ edges with $2$ vertices in $abcd$.
Furthermore, equality can only hold if deleting $abcd$ leaves $\binom{5}{3}$ edges,
i.e.\ a copy of $K^3_5$. If equality does not hold then we can find a copy of $K^3_5$
as above, so either way we have a copy of $K^3_5$. Let $X$ be the vertex set of this $K^3_5$.
Now note that for any $Y$ spanning a copy of $K^3_4$
we either have $i(Y) \ge b(n)-b(n-4)+1$ or $G \sm Y$ spans $K^3_5$.
Also, there cannot be $3$ vertices $x_1,x_2,x_3$ in $X$ such that $G \sm (X \sm x_i)$
spans $K^3_5$ for $1\le i \le 3$, as then $x_1,x_2,x_3$ together with any $3$ vertices
of $G \sm X$ spans a copy of $F_{3,3}$. Now we can modify the second
calculation above to get $3(b(n)-b(n-4)+1) + 2(b(n)-b(n-4)) \le 5(10+Z+m(n-5))$,
so $Z \ge 8(n-5)-2 = 30 > 7(n-5)$. It follows that there is some $f \in V(G) \sm \{a,b,c,d,e\}$
such that at least $8$ pairs from $abcde$ form an edge with $f$.
As above, this creates a copy of $F_{3,3}$, so we have a contradiction.
This proves the theorem. \qed

\section{Concluding remarks}

The obvious unanswered question from this paper is to characterise
the extremal examples for the problem: is it true that for $n \ge 6$,
equality can only be achieved by $B(n)$? It may be that there is a simple proof
of this statement, but if not, one might still hope to prove it for $n$
sufficiently large by the stability method. The idea would be to show
that any $F_{3,3}$-free graph $G$ on $n$ vertices with $e(G) \sim b(n)$
is `structurally close' to $B(n)$. Then one would hope to show that
any construction $B'(n)$ that is sufficiently close to $B(n)$ is
suboptimal, unless $B'(n)=B(n)$. (See \cite[Section 5]{K11} for further
discussion of this method.)

Consideration of this stability question leads us in turn to the question of
what constructions are (near) extremal for the multigraph result,
Lemma \ref{multi}. We have already seen two very different constructions
that achieve the maximum $m(n)$. However, we can rule out $M_2(n)$ by
returning to the $3$-graph world. Suppose that $G$ is an $F_{3,3}$-free
$3$-graph and $M$ is the link multigraph of a $K^3_4$ $abcd$ in $G$.
Let $J$ be the set of pairs of multiplicity at least $3$ in $M$.
We claim that $J$ is $K_t$-free, where $t = R_4(3)$ is the $4$-colour
Ramsey number for triangles. To see this, colour the pairs of multiplicity $3$
according to which colour from $abcd$ is {\em not} available, and colour
the pairs of multiplicity $4$ arbitrarily from $abcd$. If $J$ contains $K_t$
then by definition we can find a monochromatic triangle $xyz$.
Without loss of generality $d$ is the missing colour for $xy$, $xz$ and $yz$.
Then $abcxyz$ is a copy of $F_{3,3}$, contradiction. Thus $J$ is $K_t$-free.

It follows that $M_2(n)$, or indeed any `sufficiently close' construction,
cannot be realised as the link multigraph of $G$. However, we can give other
multigraph constructions that are not ruled out on these grounds.
For example, take a balanced partition of $n$ vertices into $4$ parts
$W$, $X$, $Y$, $Z$, such that all pairs within a part have multiplicity $2$,
all pairs between $W$ and $X$ or between $Y$ and $Z$ have multiplicity $4$,
and all pairs between the other four pairs of parts have multiplicity $3$.
One can check that this construction satisfies the hypotheses of Lemma \ref{multi},
and the number of edges is approximately $m(n)$. However, as far as we are aware,
it does not seem to arise as the link multigraph of a near extremal $F_{3,3}$-free $3$-graph.
The potentially large variety of multigraph constructions suggests that this may be
a difficult approach to proving stability for $F_{3,3}$, so perhaps other ideas are needed.

\medskip

\nib{Acknowledgements.}
We thank the Mathematisches Forschungsinstitut Oberwolfach for their hospitality,
and the organisers (Jeff Kahn, Angelika Steger and Benny Sudakov) for inviting
us to the meeting at which this research was conducted.

\end{document}